# General Two-Variable Functions on the Slide Rule

*István Szalkai*

University of Pannonia, Veszprém, Hungary, szalkai@almos.uni-pannon.hu

**Abstract**

We discuss the general theory of realizing two-variable functions on slide rules [3] and offer some new scales for practical use.

## 1 Introduction

All slide rules have different scales [9], many of which are for special tasks. However, all of these scales represent only one-variable functions, i.e., tables; the only two-variable functions are multiplication and division. But the mechanical construction of slide rules (two moving strips) is primarily for two varying quantities, so why are there no other two-variable functions on slide rules ?

Below we shortly summarize the general theory of two-variable functions on slide rules, give many examples, and some scales at the end of the paper (See Figure 3). High resolution scales [6],[7] were created and uploaded by myself (in 2015, using Turbo Pascal 3.0).

Two-variable function scales cannot be new, but I was unable to find any publications about them.

## 2 Two Examples

**i)** Let us write the (nonnegative) real numbers $x$ at the distance $1/x$ from the beginning points $S_1$ and $S_2$ on the two scales, where $S_1$ and $S_2$ correspond to the points on the scales which have values of infinity ($\infty$). If we slip $S_2$ of scale two to the mark $x$ on scale one (at distance $1/x$) and now look for the mark $y$ (distance $1/y$) on scale two, the mark $z$ on scale one (distance $1/z$) will satisfy the equality:

$$\frac{1}{z} = \frac{1}{x} + \frac{1}{y} \qquad (1)$$

or, in *Replus* (reciprocal of plus) form:

$$z = \frac{1}{1/x + 1/y} = \frac{xy}{x+y}$$

as shown below:

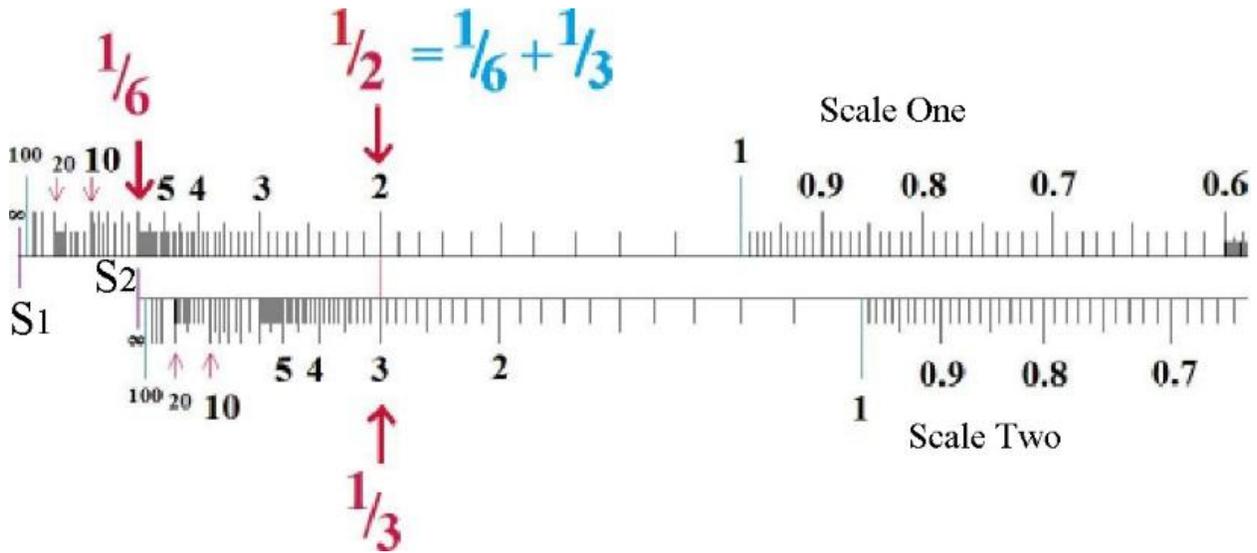

**FIGURE 1. The Reciprocal Scale**

Repeated movements clearly result in:

$$\frac{1}{z} = \frac{1}{x_1} + \cdots + \frac{1}{x_n}$$

*Replus* (reciprocal plus) is often used in optics for optical power, in electricity for parallel resistors, in parallel jobs (more workers helping each other), in geometry (incircles (inscribed circles) of triangles), and calculating harmonic means (after multiplying by 2 or by *n* for *n* numbers).

**ii)** Another famous example is *Quadplus* (quadratic plus) for *Pythagoras' theorem*.

For the formula:

$$z^2 = x^2 + y^2 \qquad (2)$$

write the (nonnegative) real numbers $x$ at the distance $x^2$ from the beginning $S_1$ and $S_2$ of the two scales, $S_1$ and $S_2$ correspond to $0^2 = 0$. If we slip $S_2$ of scale <u>two</u> to the mark $x$ on scale <u>one</u> (distance $x^2$) and now look for the mark $y$ (distance $y^2$) on scale <u>two</u>, the mark $z$ on scale <u>one</u> (distance $z^2$) will satisfy (2), or equivalently:

$$z = \sqrt{x^2 + y^2} \qquad (3)$$

in one movement:

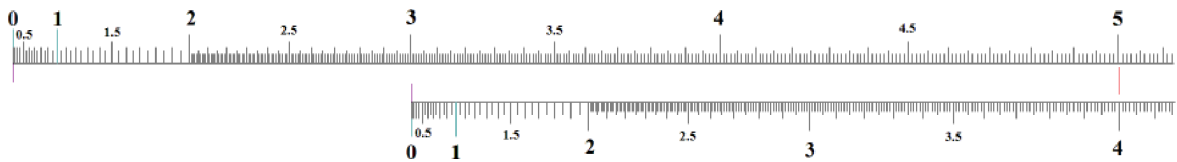

**FIGURE 2.** *The Quadratic Scale*

(2) and (3) provideus with higher <u>dimensional</u> Pythagoras' theorem and the <u>quadratic mean</u>, because several movements give us:

$$z = \sqrt{x_1^2 + \cdots + x_n^2} \qquad (4)$$

For the quadratic mean we also have to divide (4) by $\sqrt{n}$. (For a story about the quadratic mean see the Addendum [4]).

Usual slide rules cannot do addition, so *Replus* and *Quadplus* cannot be computed easily, while the above scales do these tasks in one movement!

In Figure 3 we enclose reciprocal and quadratic scales in A4 size, and we offer also the high resolution ones [7],created and uploaded by myself in 2015.

### 3 In General

Initially denote the left hand side ends of the two strips (i.e., the beginnings of their scales) by $S_1$ and $S_2$. Choose any two <u>strictly monotonic</u> functions $f$ and $g$ and write the (possible) values $x, y \in R$ on the two strips at the distances $f(x)$ and $g(y)$ from $S_1$ and $S_2$ respectively Now, sliding the second strip's $S_2$ mark to the mark $x$ on the first strip, and looking to the mark $y$ on the second strip, we are at the distance $f(x)+g(y)$ from $S_1$. If a third scale has its beginning $S_3 = S_1$ and the numbers $z \in R$ are similarly written at distance $F(z)$ then we have, in general:

$$F(z) = f(x) \pm g(y) \qquad (5)$$

or

$$z = F^{-1}(f(x) \pm g(y)) \qquad (6)$$

Clearly in the case $F=f$, in practice, we need only two scales.

Other functional relations may also be transformed into (5). For example:

$$a \cdot u(x) \cdot v(y) + b \cdot u(x) + c \cdot v(y) + d \cdot w(z) + e = 0 \qquad (7)$$

for $a,b,c,d,e \in R$, $a \neq 0$ is equivalent to:

$$\left(u(x) + \frac{c}{a}\right)\left(v(y) + \frac{b}{a}\right) = \frac{bc}{a^2} - \frac{dw(z)+e}{a} \qquad (8)$$

and to:

$$\log\left(u(x) + \frac{c}{a}\right) + \log\left(v(y) + \frac{b}{a}\right) = \log\left(\frac{bc}{a^2} - \frac{dw(z)+e}{a}\right) \qquad (9)$$

which is in the form of (5), assuming all the quantities in the brackets are positive. Obviously strictly monotonic $u(x), v(y)$, and $w(z)$ ensure that $f(x), g(y)$, and $F(z)$ in (9) are also strictly monotonic, i.e., are suitable for using with a sliderule.

Similarly, the expression

$$u(x) \cdot v(x) \cdot w(x) + u(x) + v(x) + w(x) = 0$$

is equivalent to

$$\frac{1-u(x)}{1+u(x)} \cdot \frac{1-v(x)}{1+v(x)} \cdot \frac{1-w(x)}{1+w(x)} = 1$$

(just multiply by the denominators), and so to

$$log\left(\frac{1-u(x)}{1+u(x)}\right) + log\left(\frac{1-v(x)}{1+v(x)}\right) = -log\left(\frac{1-w(x)}{1+w(x)}\right),$$

which is if form (5). Many other formulae, too can be transformed to (5) [17].

## 4 More Examples

**iii)** First of all, similarly to *Replus* (1) and *Quadplus* (2), all the scales $x^\alpha$ for each $\alpha \in R$, $\alpha \neq 0$ can be realized and make easier the calculation of $z^\alpha = x^\alpha \pm y^\alpha$ and of the <u>general mean</u>: [10],[13] .

$$z = H_\alpha = \sqrt[\alpha]{\frac{x_1^\alpha + \cdots + x_n^\alpha}{n}} \quad (10)$$

(after dividing $\sqrt[\alpha]{x_1^\alpha + \cdots + x_n^\alpha}$ by $\sqrt[\alpha]{n}$). For example, the radii of three circles, tangential to each other and a line satisfy the equation for $\alpha = -1/2$: [1],[2],[13].

$$\frac{1}{\sqrt{r_1}} = \frac{1}{\sqrt{r_2}} + \frac{1}{\sqrt{r_3}}$$

See also Descartes' theorem [14].

$H_\alpha$ for various $\alpha$ is often used in physics (springs, inductors, resistors, optics), in probability theory (dispersion) and in geometry [13]. For another example, the inductivity of serial inductors is

$$\sqrt{L_3} = \sqrt{L_1} + \sqrt{L_2}$$

We leave the question "*How is it possible to realize all the general means* $H_\alpha$ [10],[13] ($\alpha \in R$) *in* (10) *using only a few scales*?" open to the Readers, until our next paper [5] appears (soon!).

Below we recall some examples [3] and give the corresponding functions *f,g* and *F* for (5) and (6). (The roles of *F* and $F^{-1}$ in our original paper [3] are interchanged, compared to (5) and (6).)

**iv)** Solving the equality $x^2+px+q=0$, the half-part

$$\sqrt{\frac{p^2}{4} - q}$$

of the solver can be realized by $f(p)=p^2/4$, $g(q)=q$ and $F^{-1}(z)=\sqrt{z}$, i.e., $F(z)=z^2$.

**v)** The half-part

$$\sqrt{\frac{p^3}{27} + \frac{q^2}{4}}$$

of the solver for cubic equations needs $f(p)=p^3/27$, $g(q)=q^2/4$ and $F^{-1}(z)=\sqrt{z}$, i.e., $F(z)=z^2$.

**vi)** The *Lorentz* transformation:

$$\frac{M}{\sqrt{1-\frac{v^2}{c^2}}}$$

is realized by $f(M)=log(M)$, $g(v)=\frac{-1}{2}log\left(1-\frac{v^2}{c^2}\right)$ and $F(z)=log(z)$.

**vii)** $n!/k!$ and $n! \cdot k!$ can be computed by $f(x)=g(x)=log(x!)$ (for non-integer $x$ we mean $log(\Gamma(x+1))$ [11], [12] and $F(z)=log(z)$. However, I am still unable to find scales for calculating $\binom{n}{k}$ in one movement, so this question is still unsolved, except [18].

**viii)** The functions
$$f(t) = R \cdot arccos\left(\frac{R}{R+t}\right),$$
$$g(h) = R \cdot arccos\left(\frac{R}{R+h}\right)$$
and
$$F(z) = z$$

help us to compute the following everyday problem: if the observer (sailor) is at height *h* (on the mast) observes the object (top of ther tower) of height *t* , then their distance is [15]  $F(z)=f(h)+g(t)$.

**ix)** The two-variable function

$$z=x^y$$

has well-known transformations to  $log(z) = y \cdot log(x)$
and
$$log(log(z)) = log(y) + log(log(x)) ,$$

so we need scales  $F(z)=log(log(z))$, $f(x)=log(log(x))$ and  $g(y)=log(y)$.  These scales can be found on most advanced (not base) slide rules: the $log(log())$ scales are called "$e^x$" or LL3 scales, while the $log()$ scale is the base scale: called "*x*" or  C or D scale.

**x)**  $T(S,R) = (aS+b)^\alpha \cdot (cR+d)^\beta$  [8],
$z = x^\alpha y^\beta$  (many equations in physics),
$z = log_a(b)$,
$V = r^2 \pi m/3$,
*Law of Sines*(*Snellius-Descartes* law),

Formulas of the form (7) and many other two- or more variable formulas. Any other ideas and formulas from the Readers are welcome!

## 5 Addendum

Finally, let us tell a story about joiners, carpenters, and an application of the quadratic mean [4], Any wooden frame is made from slats *a,a,b,b*. Cutting precise lengths is easy, but the right angle is harder. So,they measure the actual diagonals *e,f* of the (plain) parallelogram, and they calculate the diagonal *d* of the desired rectangle as *d*=(*e*+*f*)/2, the arithmetical mean. However, two easy applications of the law of cosines show that actually *d* is the quadratic mean $d = \sqrt{\frac{e^2+f^2}{2}}$. That is why I think a slide rule with the scale shown in Figure 2 would help them!

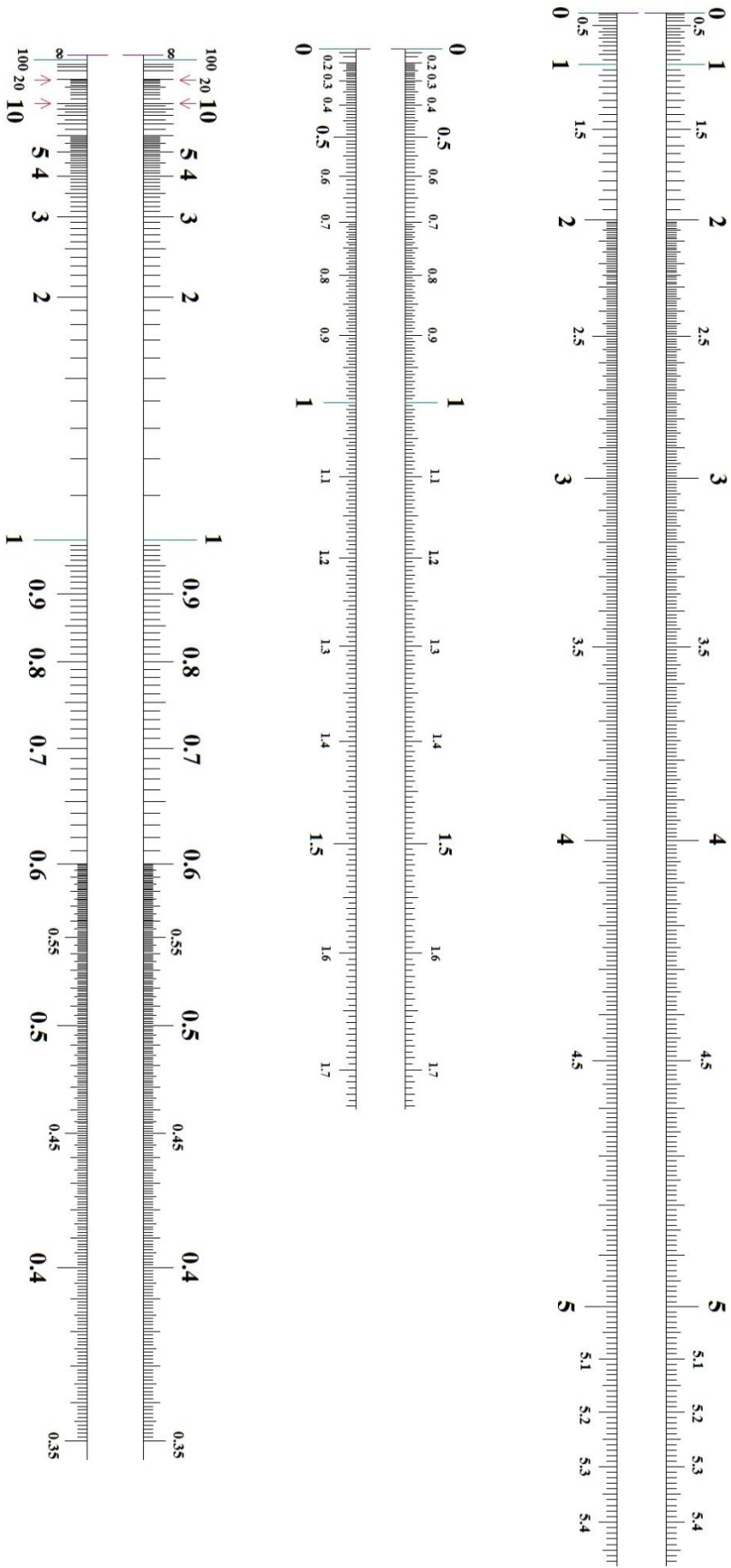

**FIGURE 3. Full Reciprocal and Quadratic Scales**